\documentclass[12pt]{article}
\usepackage{amsfonts}
\usepackage{graphicx}
\usepackage{amsmath}
\usepackage{amssymb}
\usepackage{psfrag}
\newcommand{\CF}{{\cal F}}

\newcommand{\CR}{\mathbb{R}}
\newcommand{\CE}{\mathbb{E}}
\newcommand{\CV}{\mathbb{V}}
\newcommand{\Ce}{{\cal E}}

\title{General laws of large numbers under sublinear expectations\thanks{This work has been supported in part by the National Basic
Research Program of China (973 Program) (Grant No. 2007CB814901 )
(Financial Risk). First version: Oct. 20, 2010. This is the third
version.}}

\author{ Feng  Hu\thanks{\small{\it E-mail address}:
hufengqf\symbol{64}163.com (F. Hu).}\\
School of Mathematics
\\ Shandong University\\Jinan 250100, China}

\date{}
\begin{document}
\maketitle

\begin{center}{\bf Abstract}\end{center}
\begin{center}
\parbox{0.8\hsize}{ \parindent=0.5 cm
In this paper, under some weaker conditions, we give three laws of
large numbers under sublinear expectations (capacities), which
extend Peng's law of large numbers under sublinear expectations in
[8] and Chen's strong law of large numbers for capacities in [1].
It turns out that these theorems are natural extensions of the
classical strong (weak) laws of large numbers to the case where
probability measures are no longer additive.}
\end{center}

\section{ \textbf {Introduction}}

  The classical strong (weak) laws of large numbers (strong (weak) LLN) as
fundamental limit theorems in probability theory play a fruitful
role in the development of probability theory and its
applications. However, these kinds of limit theorems have always
considered additive probabilities and additive expectations. In
fact, the additivity of probabilities and expectations has been
abandoned in some areas because many uncertain phenomena can not
be well modelled by using additive probabilities and additive
expectations. Motivated by some problems in mathematical
economics, statistics and quantum mechanics, a number of papers
have used non-additive probabilities (called capacities) and
nonlinear expectations (for example Choquet integral/expectation,
$g$-expectation) to describe and interpret the phenomena.
Recently, motivated by the risk measures, super-hedge pricing and
model uncertainty in finance, Peng [5-9] initiated the notion of
independent and identically distributed (IID) random variables
under sublinear expectations. Furthermore, he showed law of large
numbers (LLN) and central limit theorem (CLT) under sublinear
expectations. In [1], Chen presented a strong law of large numbers
for capacities induced by sublinear expectations with the notion
of IID random variables initiated by Peng.

The purpose of this paper is to investigate one of the very
important fundamental results in the theory of Peng's sublinear
expectations: tlaw of large numbers. All of the results on laws of
large numbers in [1] and [8] require that the sequence of ¡°random
variables¡± is independent and identically distributed. In this
paper we intend to obtain three laws of large numbers without the
requirement of identical distribution. Under some weaker
conditions, we prove three laws of large numbers under Peng's
sublinear expectations, which extend Peng's law of large numbers
under sublinear expectations in [8] and Chen's strong law of large
numbers for capacities in [1].

This paper is organized as follows: in Section 2, we recall some
notions and lemmas under sublinear expectations. In Section 3, we
give our main results including the proofs.

\section{Notions and Lemmas}

In this section, we present some preliminaries in the theory of
sublinear expectations.\\

\noindent{\bf DEFINITION 2.1} (see [5-9]). Let $\Omega$ be a given
set and let ${\cal H}$ be a linear space of real valued functions
defined on $\Omega$. We assume that all constants are in ${\cal
H}$ and that $X\in{\cal H}$ implies $|X|\in{\cal H}$. ${\cal H}$
is considered as the space of our ''random variables''. A
nonlinear expectation $\mathbb{E}$ on ${\cal H}$ is a functional
$\mathbb{E}$ : ${\cal H}\mapsto\mathbb R$ satisfying the following
properties: for all $X$, $Y\in{\cal H}$, we have\\
(a) Monotonicity:\ \ \ \ If $X\geq Y$ then
$\mathbb{E}[X]\geq\mathbb{E}Y]$.\\
(b) Constant preserving: $\mathbb{E}[c]=c$.\\
The triple $(\Omega,{\cal H},\mathbb{E})$ is called a nonlinear
expectation space (compare with a probability space $(\Omega,{\cal
F},P)$). We are mainly concerned with sublinear expectation where
the expectation $\mathbb{E}$ satisfies
also\\
(c) Sub-additivity:\ \ \ \ $\mathbb{E}[X]-\mathbb{E}[Y]\leq
\mathbb{E}[X-Y]$.\\
(d)Positive homogeneity: $\mathbb{E}[\lambda
X]=\lambda\mathbb{E}[X]$, $\forall\lambda\geq0$.\\
If only (c) and (d) are satisfied, $\mathbb{E}$ is called a
sublinear functional.\\

The following representation theorem for sublinear expectations is
very useful (see Peng [8, 9] for the proof).\\

\noindent{\bf LEMMA 2.1.} Let $\mathbb{E}$ be a sublinear functional
defined on $(\Omega,{\cal H})$, i.e., (c) and (d) hold for
$\mathbb{E}$. Then there exists a family
$\{E_\theta:\theta\in\Theta\}$ of linear functionals on
$(\Omega,{\cal H})$ such that
$$\mathbb{E}[X]=\max\limits_{\theta\in\Theta}E_\theta[X].$$ If (a) and
(b) also hold, then $E_\theta$ are linear expectations for
$\theta\in\Theta$. If we make furthermore the following assumption:
(H) For each sequence $\{X_n\}_{n=1}^\infty\subset{\cal H}$ such
that $X_n(\omega)\downarrow0$ for $\omega$, we have
$\hat{E}[X_n]\downarrow0$. Then for each $\theta\in\Theta$, there
exists a unique ($\sigma$-additive) probability measure $P_\theta$
defined on $(\Omega,\sigma({\cal H}))$ such that
$$E_\theta[X]=\int_\Omega X(\omega){\rm d}P_\theta(\omega),\ \ X\in{\cal
H}.$$\\

\noindent{\bf REMARK 2.1.} Lemma 2.1 shows that under $(H)$,
indeed, a stronger representation holds.  That is, if $\mathbb{E}$
is  a sublinear expectation on ${\cal H}$ satisfying (H), then
there exists a set (say $\hat{\cal P}$) of probability measures
such that
$$\mathbb{E}[\xi]=\sup_{P\in\hat{\cal P}}E_P[\xi],\ \ \
-\mathbb{E}[-\xi]=\inf_{P\in\hat{\cal P}} E_{P}[\xi].$$ Therefore,
without confusion, we sometimes call supremum expectations as
sublinear expectations.\\

Given a sublinear expectation $\CE$, let us denote the conjugate
expectation $\Ce$ of  sublinear $\CE$  by
$$
\Ce[X]:=-\CE[-X], \quad  \forall X\in \mathcal {H}
$$
Obviously, for all $X\in\mathcal {H},$ $\Ce[X]\leq \CE[X].$

Furthermore, let us denote a pair $(\CV,v)$ of capacities by
$$
\CV(A):=\CE[I_A], \quad  v(A):= \Ce[I_A],\quad \forall A\in\CF.
$$
It is easy to check that
$$
\CV(A)+v(A^c)=1,\quad \forall A\in\CF
$$
where $A^c$ is the complement set of $A.$\\

\noindent{\bf DEFINITION 2.2.} A set function $V$:
$\CF\rightarrow[0,1]$ is called upper continuous capacity if it
satisfies
\begin{description}
\item {(1)} $V(\phi)=0,V(\Omega)=1$. \item {(2)} $V(A)\leq V(B),$
whenever $A\subset B$ and $A,B\in\CF$. \item {(3)}
$V(A_{n})\downarrow V(A)$, if $A_{n}\downarrow A$, where $A_{n},
A\in \CF$.
\end{description}

\noindent{\bf ASSUMPTION A.} Throughout this paper, we assume that
$\CE$ is a sublinear expectation, $\CV$ is an upper continuous
capacity generated
by sublinear expectation $\CE$.\\

The following is the notion of IID random variables under sublinear
expectations  introduced by Peng [5-9].\\

\noindent{\bf  DEFINITION 2.3.} {\bf Independence:} Suppose that
$Y_1,Y_2,\cdots,Y_n$ is a sequence of random variables such that
$Y_i \in\mathcal {H}.$ Random variable $Y_n$ is said to be
independent to $X:=(Y_1,\cdots,Y_{n-1})$ under $\mathbb{E}$, if for
each measurable function $\varphi$ on $R^n$ with $\varphi(X,Y_n)\in
\mathcal {H}$ and $\varphi(x,Y_n)\in\mathcal {H}$ for each $x\in
\mathbb{R}^{n-1},$ we have
$$\mathbb{E}[\varphi(X,Y_n)]=\mathbb{E}[\overline{\varphi}(X)],$$
where $\overline{\varphi}(x):=\mathbb{E}[\varphi(x,Y_n)]$ and
$\overline{\varphi}(X)\in\mathcal {H}$.

{\bf Identical distribution:} Random variables $X$ and $Y$ are said
to be identically distributed, denoted by $X\overset{d}{=}Y$, if for
each $\varphi$ such that $\varphi(X), \; \varphi(Y)\in \mathcal
{H}$,
$$\mathbb{E}[\varphi(X)]=\mathbb{E}[\varphi(Y)].$$

{\bf Sequence of IID random variables:} A sequence of IID random
sequence $\{X_i\}_{i=1}^\infty$ is said to be  IID random variables,
if $X_i\overset{d}{=}X_1$ and $X_{i+1}$ is independent to
$Y:=(X_1,\cdots,X_i)$ for each $i\ge 1.$\\

The following lemma shows the relation between Peng's independence
and pairwise independence in Maccheroni and Marinacci [3] and Marinacci [4].\\

\noindent{\bf LEMMA 2.2} (see Chen [1]). Suppose that $X,Y \in
\mathcal {H}$ are two random variables. $\CE$ is a sub-linear
expectation and $(\CV,v)$ is the pair of capacities generated by
$\CE.$ If  random variable $X$ is independent to $Y$ under $\CE$,
then $X$ also is pairwise independent to $Y$ under capacities $\CV,$
and $v$ e.g. for all subsets $D$ and $G\in{\cal B}(\CR),$
$$
V(X\in D, Y\in G)=V(X\in D)V(Y\in G)
$$
holds for both capacities $\CV$ and $v$.\\

Borel-Cantelli Lemma is still true for capacity under some
assumptions.\\

\noindent{\bf LEMMA 2.3} (see Chen [1]). Let $\{A_n,n\geq1\}$ be a
sequence of events in $\CF$ and $(\CV,v)$ be a pair of capacities
generated by sublinear expectation $\CE$.
 \begin{description}
 \item{(1)} If $\sum\limits_{n=1}^\infty\CV(A_n)<\infty,$ then $\CV\left(\bigcap\limits_{n=1}^\infty\bigcup\limits_{i=n}^\infty
 A_i\right)=0.$
\item{(2)} Suppose that  $\{A_n,n\geq1\}$ are pairwise independent
with respect to $v$, e.g. $$v\left(\bigcap\limits_{i=1}^\infty
A_i^c\right)=\prod_{i=1}^\infty v\left(A_i^c\right).$$ If
$\sum\limits_{n=1}^\infty{\CV}(A_n)=\infty $, then
$$\CV\left(\bigcap\limits_{n=1}^\infty\bigcup\limits_{i=n}^\infty A_i\right)=1.$$
\end{description}

\noindent{\bf  DEFINITION  2.4 ( Maximal distribution)} (see Peng
[8, 9]). Let $C_{b,Lip}(\mathbb R)$ denote the space of bounded
and Lipschitz continuous functions. A random variable $\eta$ on
sublinear expectation space $(\Omega,{\cal H},\mathbb E)$ is
called maximal distributed if
$$\mathbb E[\varphi(\eta)]=\sup\limits_{\underline\mu\leq y\leq\overline\mu}\varphi(y),
\ \ \ \forall\varphi\in C_{b,Lip}(\mathbb R),$$ where
$\overline\mu:=\mathbb E[\eta]$ and $\underline\mu:=\Ce[\eta]$.\\

\noindent{\bf REMARK  2.2} (see Peng [8, 9]). Let $\eta$ be
maximal distributed with $\overline{\mu}=\CE[\eta]$,
$\underline{\mu}=\Ce[\eta]$, the distribution of $\eta$ is
characterized by the following parabolic PDE: $$\partial_t u
-g(\partial_{x} u)=0,\quad  u(0,x)=\varphi(x),$$ where
$u(t,x):=\CE[\varphi(x+t\eta)],$ $(t,x)\in[0,\infty)\times \CR$,
$g(x):=\overline\mu x^+-\underline\mu x^-$ and $x^+:=\max\{x,0\}$,
$x^-:=(-x)^+$.\\

With the notion of IID under sublinear expectations, Peng shows a
law of large numbers under sublinear expectations (see Theorem 5.1
in
Peng [8]).\\

\noindent{\bf LEMMA 2.4 (Law of large numbers under sublinear
expectations)}. Let $\{X_i\}_{i=1}^\infty$ be a sequence of IID
random variables with finite means $\overline\mu=\CE[X_1], \;\;
\underline \mu=\Ce[X_1].$  Suppose $\CE[|X_1|^{2}] <\infty$. Then
for any continuous and linear growth function $\varphi,$
$$
\CE\left[\varphi\left(\frac 1n \sum_{i=1}^nX_i\right)\right]\to
\sup_{\underline \mu\leq y\leq\overline \mu}\varphi(y), \;  \hbox{
as}\; n\to \infty.$$\\

The following lemma is a strong law of large numbers for
capacities, which can be found in Chen [1].\\

\noindent{\bf LEMMA 2.5 ( Strong law of large numbers for
capacities)}. Let $\{X_i\}_{i=1}^\infty$ be a sequence of IID
random variables for sublinear expectation $\CE$. Suppose
$\CE[|X_1|^{2}]<\infty$. Set $\overline\mu:=\CE[X_1],$ $\underline
\mu:=\Ce[X_1]$  and $ S_n:=\sum\limits_{i=1}^n X_i.$ Then
\begin{description}
\item{(I)}
\begin{equation*}
v\left( \underline \mu\leq \liminf\limits_{n\to\infty} S_n/n \leq
\limsup\limits_{n\to\infty} S_n/n \leq \overline \mu\right)=1.
\end{equation*}
\item{(II)}
$$\CV\left( \limsup\limits_{n\to\infty}S_n/n= \overline
\mu\right)=1,\quad \CV\left( \liminf\limits_{n\to\infty} S_n/n=
\underline \mu\right)=1.
$$
\item{(III)} $\forall b\in[\underline\mu,\overline\mu]$,
 $$
 \CV\left (\liminf\limits_{n\to\infty}|S_n/n-b|=0\right)=1.
 $$
\end{description}

\section{Main Results}

\noindent{\large\bf 3.1 General law of large numbers under
sublinear
expectations}\\

\noindent{\bf THEOREM 3.1 (General law of large numbers under
sublinear expectations)}. Let a sequence $\{X_i\}_{i=1}^\infty$
which is in a sublinear expectation space $(\Omega,{\cal
H},\mathbb E)$ satisfy the following conditions:

\noindent(\rm i) each $X_{i+1}$ is independent of $(X_1,\cdots,
X_i)$, for $i = 1,2,\cdots$;

\noindent(\rm ii) $\mathbb E[X_i]=\overline{\mu}_i$,
$\Ce[X_i]=\underline{\mu_i}$, where
$-\infty<\underline{\mu_i}\leq\overline{\mu}_i<\infty$;

\noindent(\rm iii) there are two constants ${\overline{\mu}}$ and
${\underline{\mu}}$ such that
$$\lim\limits_{n\rightarrow\infty}\frac{1}{n}\sum\limits_{i=1}^n|\underline{\mu_i}-\underline{\mu}|=0,
\ \ \
\lim\limits_{n\rightarrow\infty}\frac{1}{n}\sum\limits_{i=1}^n|\overline{\mu_i}-\overline{\mu}|=0;$$

\noindent(iv) $\sup\limits_{i\ge 1}\CE[|X_i|^{2}]<\infty$. Then for
any continuous and linear growth function $\varphi,$
$$\CE\left[\varphi\left(\frac 1n \sum_{i=1}^nX_i\right)\right]\to
\sup_{\underline \mu\leq y\leq\overline\mu}\varphi(y), \;  \hbox{
as}\; n\to \infty.$$

\noindent{\bf PROOF.} The main idea comes from Theorem 3.1 Li and
Shi [2]. First we prove the case that $\varphi$ is a bounded and
Lipschitz continuous function. For a small but fixed $h>0$, Let
$V$ be the unique viscosity solution of the following equation
\begin{equation}\label{ED}
\partial_t V +g(\partial_{x}V)=0,\ \ \ (t,x)\in[0,1+h]\times \CR,\ \ \
V(1+h,x)=\varphi(x),\end{equation} where $g(x):=\overline\mu
x^+-\underline\mu x^-$. According to the definition of maximal
distribution, we have
$$V(t,x)=\CE[\varphi(x+(1+h-t)\eta)].$$ Particularly,
\begin{equation}\label{ED}V(h,0)=\CE[\varphi(\eta)],\ \ \ V(1+h,x)=\varphi(x).\end{equation}
Since (1) is a uniformly parabolic PDE, by the interior regularity
of $V$ (see Wang [10]), we have
$$||V||_{C^{1+\alpha/2,1+\alpha}([0,1]\times\CR)}<\infty,\ \ \ \hbox{for some}\ \ \ \alpha\in(0,1).$$

We set $\delta:=\frac{1}{n}$ and ${S}_0:=0$. Then
$$\begin{array}{lcl}&&V(1,\delta{S}_n)-V(0,0)=\sum\limits_{i=0}^{n-1}
\{V((i+1)\delta,\delta{S}_{i+1})-V(i\delta,\delta{S}_{i})\}\\
&&=\sum\limits_{i=0}^{n-1}\left\{\left[V((i+1)\delta,\delta{S}_{i+1})-V(i\delta,\delta{S}_{i+1})\right]
+\left[V(i\delta,\delta{S}_{i+1})-V(i\delta,\delta{S}_{i})\right]\right\}\\
&&=\sum\limits_{i=0}^{n-1}\{I_\delta^i+J_\delta^i\},
\end{array}$$ with, by Taylor's expansion, $$J_\delta^i=\partial_t V(i\delta,\delta{S}_i)\delta
+\partial_x V(i\delta,\delta{S}_i)X_{i+1}\delta,$$

$$\begin{array}{lcl}I_\delta^i&=&\int_0^1\left[\partial_t
V((i+\beta)\delta,\delta{S}_{i+1}) -\partial_t
V(i\delta,\delta{S}_{i+1})\right]{\rm
d}\beta\delta\\
&+&\left[\partial_t V(i\delta,\delta{S}_{i+1})-\partial_t
V(i\delta,\delta{S}_{i})\right]\delta\\
&+&\int_0^1\left[\partial_{x} V(i\delta,\delta{S}_{i}+\beta\delta
X_{i+1})-\partial_{x} V(i\delta,\delta{S}_{i})\right]{\rm d}\beta
X_{i+1}\delta.\end{array}$$ Thus
\begin{equation}\label{ED}\CE\left[\sum\limits_{i=0}^{n-1}J_\delta^i\right]+\Ce\left[\sum\limits_{i=0}^{n-1}I_\delta^i\right]\leq\CE[V(1,\delta{S}_n)]-V(0,0)\leq
\CE\left[\sum\limits_{i=0}^{n-1}J_\delta^i\right]+\CE\left[\sum\limits_{i=0}^{n-1}I_\delta^i\right].\end{equation}
From (1) as well as the independence of $X_{i+1}$ to $(X_1,\cdots,
X_i)$, it follows that

$$\begin{array}{lcl}&&\CE[J_\delta^i]=\CE\left[\partial_t V(i\delta,\delta{S}_i)\delta
+\partial_x
V(i\delta,\delta{S}_i)X_{i+1}\delta\right]\\
&&=\CE\left\{\partial_t
V(i\delta,\delta{S}_i)\delta+\delta[(\partial_xV(i\delta,\delta{S}_i))^{+}\overline{\mu_{i+1}}-
(\partial_xV(i\delta,\delta{S}_i))^{-}\underline{\mu_{i+1}}]\right\}\\
&&\leq\CE\left\{\partial_t
V(i\delta,\delta{S}_i)\delta+\delta[(\partial_xV(i\delta,\delta{S}_i))^{+}\overline{\mu}-
(\partial_xV(i\delta,\delta{S}_i))^{-}\underline{\mu}]\right\}\\
&&+\delta\CE[(\partial_xV(i\delta,\delta{S}_i))^{+}(\overline{\mu_{i+1}}-\overline{\mu})+
(\partial_xV(i\delta,\delta{S}_i))^{-}(\underline{\mu_{i+1}}-\underline{\mu})]\\
&&=\delta\CE[(\partial_xV(i\delta,\delta{S}_i))^{+}(\overline{\mu_{i+1}}-\overline{\mu})+
(\partial_xV(i\delta,\delta{S}_i))^{-}(\underline{\mu_{i+1}}-\underline{\mu})]\\
&&\leq\delta(|\overline{\mu_{i+1}}-\overline{\mu}|+|\underline{\mu_{i+1}}-\underline{\mu}|)\CE[|\partial_xV(i\delta,\delta{S}_i)|].\end{array}$$
But since both $\partial_t V$ and $\partial_{x} V$ are uniformly
$\alpha$-h\"{o}lder continuous in $x$ and
$\frac{\alpha}{2}$-h\"{o}lder continuous in $t$ on
$[0,1]\times\CR$, it follows that
$$|\partial_xV(i\delta,\delta{S}_i)-\partial_xV(0,0)|\leq
C|\delta{S}_i|^\alpha+|i\delta|^\frac{\alpha}{2},$$ where $C$ is
some positive constant. Since
$$\CE[|\delta{S}_i|^\alpha]\leq\CE[|\delta{S}_i|]+1\leq\sup\limits_{i\ge 1}\CE[|X_i|]+1.$$
Hence, by (\rm iv), we claim that there exists a constant $C_1 > 0$,
such that
$$\CE[|\partial_xV(i\delta,\delta{S}_i)|]\leq C_1.$$ Then we obtain
$$\CE\left[\sum\limits_{i=0}^{n-1}J_\delta^i\right]\leq\sum\limits_{i=0}^{n-1}
\CE\left[J_\delta^i\right]\leq
C_1\frac{1}{n}\sum\limits_{i=0}^{n-1}(|\overline{\mu_{i+1}}-\overline{\mu}|+|\underline{\mu_{i+1}}-\underline{\mu}|).$$
From (\rm iii), we have
$\limsup\limits_{n\rightarrow\infty}\CE\left[\sum\limits_{i=0}^{n-1}J_\delta^i\right]\leq0.$

In a similar manner of the above, we also have
$$\CE\left[\sum\limits_{i=0}^{n-1}J_\delta^i\right]\geq-\sum\limits_{i=0}^{n-1}\delta
(|\overline{\mu_{i+1}}-\overline{\mu}|+|\underline{\mu_{i+1}}-\underline{\mu}|)\CE[|\partial_xV(i\delta,\delta{S}_i)|].$$
By (\rm iii), it follows that
$\liminf\limits_{n\rightarrow\infty}\CE\left[\sum\limits_{i=0}^{n-1}J_\delta^i\right]\geq0.$
So we can claim that
\begin{equation}\label{ED}\lim\limits_{n\rightarrow\infty}\CE\left[\sum\limits_{i=0}^{n-1}J_\delta^i\right]=0.\end{equation}

For $I_\delta^i$, since both $\partial_t V$ and $\partial_{x} V$ are
uniformly $\alpha$-h\"{o}lder continuous in $x$ and
$\frac{\alpha}{2}$-h\"{o}lder continuous in $t$ on $[0,1]\times\CR$,
we then have $$|I_\delta^i|\leq
C\delta^{1+\alpha/2}(1+|X_{i+1}|^\alpha+|X_{i+1}|^{1+\alpha}).$$ It
follows that
$$\CE[|I_\delta^i|]\leq C\delta^{1+\alpha/2}(1+\CE[|X_{i+1}|^\alpha]+\CE[|X_{i+1}|^{1+\alpha}]).$$
Thus
\begin{equation}\label{ED}\begin{array}{lcl}&&-C(\frac{1}{n})^{\frac{\alpha}{2}}\frac{1}{n}\sum\limits_{i=0}^{n-1}\left(1+
\CE\left[|X_{i+1}|^\alpha\right]+\CE\left[|X_{i+1}|^{1+\alpha}\right]\right)+\CE\left[\sum\limits_{i=0}^{n-1}J_\delta^i\right]\\
&&\leq \CE[V(1,\delta{S}_n)]-V(0,0)\\
&&\leq
C(\frac{1}{n})^{\frac{\alpha}{2}}\frac{1}{n}\sum\limits_{i=0}^{n-1}\left(1+
\CE\left[|X_{i+1}|^\alpha\right]+\CE\left[|X_{i+1}|^{1+\alpha}\right]\right)+\CE\left[\sum\limits_{i=0}^{n-1}J_\delta^i\right].\end{array}\end{equation}
Therefore, from (\rm iv) and (4), as $n\rightarrow\infty$, we have
\begin{equation}\label{ED}\lim\limits_{n\rightarrow\infty}\CE[V(1,\delta S_n)]=
V(0,0).\end{equation}

On the other hand, for each $t$, $t^{'}\in[0,1+h]$ and $x\in\CR$,
$$|V(t,x)-V(t^{'},x)|\leq C|t-t^{'}|.$$ Thus
\begin{equation}\label{ED} |V(0,0)-V(h,0)|\leq Ch\end{equation} and, by
(2)
\begin{equation}\label{ED}\begin{array}{lcl}&&|\CE[V(1,\delta{S}_n)]-\CE[\varphi(\delta S_n)]|\\
&&=|\CE[V(1,{\delta}{S}_n)]-\CE[V(1+h,{\delta}S_n)|\leq
Ch.\end{array}\end{equation} It follows from (6)-(8) that
$$\limsup\limits_{n\rightarrow\infty}\left|\CE\left[\varphi\left(\frac{S_n}{{n}}\right)\right]-\CE[\varphi(\eta)]\right|
\leq2C{h}.$$ Since $h$ can be arbitrarily small, we have
$$\lim\limits_{n\rightarrow\infty}\CE\left[\varphi\left(\frac{S_n}{{n}}\right)\right]
=\CE\left[\varphi\left(\eta\right)\right]=\sup_{\underline \mu\leq
y\leq\overline\mu}\varphi(y).$$

The rest proof of Theorem 3.1 is very similar to that of Lemma 5.5
in peng [8]. So we omit it.\\

From Theorem 3.1, we easily claim the following corollary.\\

\noindent{\bf COROLLARY 3.1.} Let a sequence $\{X_i\}_{i=1}^\infty$
which is in a sublinear expectation space $(\Omega,{\cal H},\mathbb
E)$ satisfy the following conditions:

\noindent(\rm i) each $X_{i+1}$ is independent of $(X_1,\cdots,
X_i)$, for $i = 1,2,\cdots$;

\noindent(\rm ii) $\mathbb E[X_i]=\overline{\mu}_i$,
$\Ce[X_i]=\underline{\mu}_i$, where
$-\infty<\underline{\mu}_i\leq\overline{\mu}_i<\infty$;

\noindent(\rm iii) there are two constants ${\overline{\mu}}$ and
${\underline{\mu}}$ such that
$\lim\limits_{i\rightarrow\infty}\underline{\mu}_i=\underline{\mu}$,
$\lim\limits_{i\rightarrow\infty}\overline{\mu}_i=\overline{\mu};$

\noindent(iv) $\sup\limits_{i\ge 1}\CE[|X_i|^{2}]<\infty$. Then for
any continuous and linear growth function $\varphi,$
$$\CE\left[\varphi\left(\frac 1n \sum_{i=1}^nX_i\right)\right]\to
\sup_{\underline \mu\leq y\leq\overline\mu}\varphi(y), \;  \hbox{
as}\; n\to \infty.$$

\noindent{\large\bf 3.2 General strong law of large numbers for
capacities induced
by sublinear expectations}\\

\noindent{\bf THEOREM 3.2 (General strong law of large numbers for
capacities)}. Under the conditions of Theorem 3.1, then
\begin{description}
\item{(I)}
\begin{equation*}
v\left( \underline \mu\leq
\liminf\limits_{n\to\infty}\frac{\sum_{i=1}^nX_i}{n}
\leq\limsup\limits_{n\to\infty}\frac{\sum_{i=1}^nX_i}{n}\leq
\overline \mu\right)=1.
\end{equation*}
\item{(II)} $\forall b\in[\underline\mu,\overline\mu]$,
 $$
 \CV\left(\liminf\limits_{n\to\infty}\left|\frac{\sum_{i=1}^nX_i}{n}-b\right|=0\right)=1.
 $$
\item{(III)}
$$\CV\left( \limsup\limits_{n\to\infty}\frac{\sum_{i=1}^nX_i}{n}= \overline
\mu\right)=1,\quad \CV\left( \liminf\limits_{n\to\infty}
\frac{\sum_{i=1}^nX_i}{n}= \underline \mu\right)=1.$$
\end{description}

In order to prove Theorem 3.2, we need  the following lemma.

\noindent{\bf LEMMA 3.1} (see Chen [1]). In a sublinear expectation
space $(\Omega,{\cal H},\mathbb E)$, let $\{X_i\}_{i=1}^\infty$ be
a sequence of independent random variables such that
$\sup\limits_{i\ge 1}\CE[|X_i|^{2}]<\infty$. Suppose that there
exists a constant $C>0$ such that
$$|X_n-\CE[X_n]|\leq C\frac{n}{\log(1+n)}, \quad n=1,2,\cdots.$$
 Then there exists a sufficiently large number $m>1$ such that
$$\sup_{n\ge 1}\CE\left[\exp\left(\frac{m\log(1+n)}{n} S_n\right)\right]<\infty.$$
Where $S_n:=\sum\limits_{i=1}^n[X_i-\CE[X_i]].$\\

\noindent{\bf PROOF OF THEOREM 3.2.} First, it is easy to show that
(I) is equivalent to the conjunction of
\begin{equation}\label{F1}
\CV\left(\limsup\limits_{n\to\infty}\frac{\sum_{i=1}^nX_i}{n}>\overline
\mu\right)=0,
\end{equation}
\begin{equation}\label{F2}
\CV\left(\liminf\limits_{n\to\infty}\frac{\sum_{i=1}^nX_i}{n}<\underline
\mu\right)=0.
\end{equation}
Indeed, write
$A:=\left\{\limsup\limits_{n\to\infty}\frac{\sum_{i=1}^nX_i}{n}>
\overline \mu\right\},$ $B:=\left\{\liminf\limits_{n\to\infty}
\frac{\sum_{i=1}^nX_i}{n}<\underline\mu\right\},$ the equivalence
can be proved from the inequality
$$ \max\{ \CV(A),V(B)\}\leq \CV\left(A\bigcup B\right)\leq \CV(A)+\CV(B).$$

We shall turn the proofs of (\ref{F1}) and (\ref{F2}) into two
steps.

{\bf  Step 1}  Assume that there exists a constant $C>0$ such that
$|X_n-\overline\mu_n|\leq \frac{Cn}{\log(1+n)}$ for $n\ge 1.$ To
prove (\ref{F1}), we shall show that for any $\varepsilon>0,$
\begin{equation}\label{R1}
\CV\left(\bigcap_{n=1}^{\infty}\bigcup_{i=n}^{\infty}\left\{\frac{\sum_{i=1}^n(X_i-\overline{\mu_i})}{n}\ge
\varepsilon\right\}\right)=0.
\end{equation}
Indeed, by Lemma 3.1, for any $\varepsilon>0,$ let us choose $m>
1/\varepsilon$ such that
$$\sup_{n\ge 1}\CE\left[\exp\left( \frac{m\log(1+n)}{n}\sum_{i=1}^n(X_i-\overline\mu_i)\right)\right]<\infty.$$
By Chebyshev's inequality,
$$\begin{array}{lcl}
\CV\left(\frac{\sum_{i=1}^n(X_i-\overline\mu_i)}n\ge\varepsilon\right)
&=&
\CV\left(\frac{m\log(1+n)}{n}\sum\limits_{i=1}^n(X_i-\overline\mu_i)\ge
\varepsilon m\log(1+n)\right)\\
&\leq & {\rm e}^{-\varepsilon m\log(1+n)}\CE\left[\exp\left(
\frac{m\log(1+n)}{n}\sum\limits_{i=1}^n(X_i-\overline
\mu_i)\right)\right]\\
 &\leq& \frac 1{(1+n)^{\varepsilon m}}\sup\limits_{n\ge 1}\CE\left[\exp\left(
\frac{m\log(1+n)}{n}\sum\limits_{i=1}^n(X_i-\overline
\mu_i)\right)\right].
\end{array}$$
Since $\varepsilon m>1,$ and $\sup\limits_{n\ge
1}\CE\left[\exp\left(
\frac{m\log(1+n)}{n}\sum\limits_{i=1}^n(X_i-\overline
\mu_i)\right)\right]<\infty.$ It then follows from the convergence
of $\sum\limits_{n=1}^\infty \frac 1{(1+n)^{\varepsilon m}},$  we
have
$$\sum_{n=1}^\infty \CV\left(\frac{\sum_{i=1}^n(X_i-\overline\mu_i)}{n}\ge\varepsilon\right)<\infty.$$
 Using the first Borel-Cantelli Lemma, we have
$$
\CV\left(\limsup\limits_{n\to\infty}
\frac{\sum_{i=1}^n(X_i-\overline\mu_i)}{n} \geq\varepsilon\right)=0
\quad\forall\varepsilon>0,
$$
which implies
$$
\CV\left(\limsup\limits_{n\to\infty}\frac{\sum_{i=1}^nX_i}{n}>\overline
\mu\right)=0.
$$
Also
$$
v\left(\limsup\limits_{n\to\infty}\frac{\sum_{i=1}^nX_i}{n}\leq
\overline \mu\right)=1.
$$

Similarly, considering the sequence $\{-X_i\}_{i=1}^\infty.$ By step
1, it suffices to obtain
$$
\CV\left(\limsup\limits_{n\to\infty}\frac{-\sum_{i=1}^n(X_i-\underline{\mu_i})}{n}>0\right)=0.
$$
Hence,
$$
\CV\left(\liminf\limits_{n\to\infty}\frac{\sum_{i=1}^nX_i}{n}
<\underline\mu\right)=0.
$$
Also
$$
v\left(\liminf\limits_{n\to\infty}\frac{\sum_{i=1}^nX_i}{n}\ge
\underline\mu\right)=1 .$$

{\bf Step 2.} Write
$$\overline X_n:=(X_n-\overline
\mu_n)I_{\left\{|X_n-\overline\mu_n|\leq \frac
{Cn}{\log(1+n)}\right\}}-\CE\left[(X_n-\overline
\mu_n)I_{\left\{|X_n-\overline\mu_n|\leq
\frac{Cn}{\log(1+n)}\right\}}\right]+\overline\mu_n.$$ It is easy to
check that $\{\overline X_i\}_{i=1}^\infty$ satisfies the
assumptions in Lemma 3.1. Indeed, obviously for each $n\ge 1,$
$$|\overline X_n-\overline
\mu_n|\leq \frac{2Cn}{\log(1+n)}.$$

On the other hand, for each $n\ge 1,$ it easy to check that
$$|\overline X_n-\overline\mu_n|\leq |X_n-\overline\mu_n|+\CE[|X_n-\overline
\mu_n|].$$  Then, by (\rm iv),
$$\CE[|\overline X_n-\overline \mu_n|^{2}]\leq
4\left(\CE[|X_n-\overline{\mu}_n|^{2}]+(\CE[|X_n-\overline{\mu}_n|])^{2}\right)<\infty.
$$

Set $\overline S_n:=\sum\limits_{i=1}^n(\overline
X_i-\overline\mu_i),$ immediately,
\begin{equation}\label{R}
\frac{1}{n}\sum\limits_{i=1}^n(X_i-\overline\mu_i)\leq
\frac{1}{n}\overline S_n+ \frac 1n \sum_{i=1}^n|X_i-\overline
\mu_i|I_{\left\{|X_i-\overline\mu_i|>\frac{Ci}{\log(1+i)}\right\}}+\frac1n
\sum_{i=1}^n\CE\left[|X_i-\overline
\mu_i|I_{\left\{|X_i-\overline\mu_i|>\frac{Ci}{\log(1+i)}\right\}}\right].
\end{equation}
Since $\sup\limits_{i\ge 1}\CE[|X_i|^{2}]<\infty$, we have

$$\begin{array}{lcl}
\sum\limits_{i=1}^\infty\frac{\CE\left[\left|X_i-\overline
\mu_i\right|I_{\left\{|X_i-\overline\mu_i|>\frac{Ci}{\log(1+i)}\right\}}\right]}{i}
&\leq& \sum\limits_{i=1}^\infty\frac{\log(1+i)}{c
i^2} \CE[|X_i-\overline \mu_i|^{2}]\\
&\leq&C_1\sum\limits_{i=1}^\infty\frac{\log(1+i)}
{i^{2}}<\infty.\end{array}.$$ By Kronecker Lemma,
\begin{equation}\label{E1}
\frac 1n \sum_{i=1}^n\CE\left[|X_i-\overline
\mu_i|I_{\left\{|X_i-\overline\mu_i|>\frac{Ci}{\log(1+i)}\right\}}\right]\to
0.
\end{equation}
Furthermore, write
$A_i:=\left\{|X_i-\overline\mu_i|>\frac{Ci}{\log(1+i)}\right\}$ for
$i\ge 1.$ It suffices now to prove that
$$\CV\left(\bigcap_{n=1}^\infty\bigcup_{i=n}^\infty A_i\right)=0.$$
Indeed, by Chebyshev's inequality,
$$\CV\left(|X_i-\overline\mu_i|>\frac{Ci}{\log(1+i)}\right)\leq
\left(\frac{\log(1+i)}{Ci}\right)^{2}\CE[|X_i-\overline\mu_i|^{2}]
$$
Hence,
$$
\sum_{i=1}^\infty\CV\left(|X_i-\overline\mu_i|>\frac{Ci}{\log(1+i)}\right)<\infty
$$
and by the first Borel-Cantelli Lemma, we have
$\CV\left(\bigcap\limits_{n=1}^\infty\bigcup\limits_{i=n}^\infty
A_i\right)=0.$\\
This implies that $\omega\not\in
\bigcap\limits_{n=1}^\infty\bigcup\limits_{i=n}A_i,$ the sequence $$
\sum_{i=1}^n\frac{|X_i-\overline
\mu_i|I_{\left\{|X_i-\overline\mu_i|>\frac{Ci}{\log(1+i)}\right\}}}i$$
converges almost surely with respect to $\CV$ as
$n\rightarrow\infty$. Applying Kronecker Lemma again,
\begin{equation}\label{E2}
\frac 1n \sum_{i=1}^n\left(|X_i-\overline
\mu_i|I_{\left\{|X_i-\overline\mu_i|>\frac{Ci}{\log(1+i)}\right\}}\right)\to
0, \hbox{ a.s} \quad \CV.
\end{equation}

Set $\lim\sup$ on both side of (\ref{R}),  then by (\ref{E1}) and
(\ref{E2}), we have
 $$\limsup\limits_{n\to\infty}\frac{\sum_{i=1}^n(X_i-\overline\mu_i)}{n}
\leq\limsup\limits_{n\to\infty}\frac{\overline S_n}{n},\quad
\hbox{a.s.}\quad \CV.$$ Since $\{\overline X_n\}_{n=1}^\infty$
satisfies the assumption of Step 1, by Step 1,
$$\CV\left(\limsup\limits_{n\to\infty}\frac{\sum_{i=1}^n\overline
X_i}{n}>\overline \mu\right)=0.$$ Also
$$v\left(\limsup\limits_{n\to\infty}\frac{\sum_{i=1}^n\overline
X_i}{n}\leq \overline \mu\right)=1.$$

In a similar manner, we can prove $$
v\left(\liminf\limits_{n\to\infty}\frac{\sum_{i=1}^nX_i}{n}\ge
\underline\mu\right)=1 .$$ Therefore, the proof of (I) is complete.

To prove (II). Denote $S_n:=\frac{1}{n}\sum_{i=1}^nX_i$. If
$\overline{\mu}=\underline{\mu}$, it is trivial. Suppose
$\overline{\mu}>\underline{\mu},$   we prove that, for any $b\in
(\underline \mu, \overline \mu)$
$$
\CV\left(\liminf_{n\to\infty}|S_{n}/{n}-b|=0\right)=1.
$$
To do so, we only need to prove that there exists an increasing
subsequence $\{n_k\}$ of $\{n\}$ such that for any $b\in
(\underline \mu, \overline \mu)$
\begin{equation}\label{R3}
\CV\left(\bigcap_{m=1}^{\infty}\bigcup_{k=m}^{\infty}\{|S_{n_k}/{n_k}-b|\leq\varepsilon\}\right)=1.
\end{equation}
Indeed, for any $0<\varepsilon\leq \min\{ \overline\mu -b,
b-\underline \mu\},$ let us choose $n_k=k^k$ for $k\ge 1.$

Set $\overline S_n:=\sum\limits_{i=1}^n(X_i-b),$ then
$$\begin{array}{lcl}
\CV\left(\left|\frac{ S_{n_k}-S_{n_{k-1}} }
{n_k-n_{k-1}}-b\right|\leq \epsilon\right) &=& \CV\left(\left|\frac{
S_{n_k-n_{k-1}} }
{n_k-n_{k-1}}-b\right|\leq\varepsilon\right)\\
&=&\CV\left(\left|\frac{ S_{n_k-n_{k-1}}-(n_k-n_{k-1})b}
{n_k-n_{k-1}}\right|\leq\varepsilon\right) \\
&=& \CV\left(\left|\frac{\overline S_{n_k-n_{k-1}}}
{n_k-n_{k-1}}\right|\leq
\varepsilon\right)\\
&\ge& \CE\left[\phi\left(\frac{\overline S_{n_k-n_{k-1}}}
{n_k-n_{k-1}}\right)\right]
\end{array}
$$
where $\phi(x)$ is defined by
$$
\phi(x)=\left\{
\begin{array}{l}
1-e^{|x|-\varepsilon},\quad\quad |x|\leq\varepsilon;\\
0,      \quad\quad\quad\quad  \quad\;\;  |x|>\varepsilon.
\end{array}
\right.
$$
Considering the sequence of random variables
$\{X_i-b\}_{i=1}^\infty.$
 Obviously $$\CE[X_i-b]=\overline\mu_i-b,\quad \Ce[X_i-
b]=\underline{\mu_i}-b.$$ Applying Theorem 3.1, we have, as $k\to
\infty$,
$$\CE\left[\phi\left(\frac{\overline S_{n_k-n_{k-1}}} {n_k-n_{k-1}}\right)\right]\to
\sup_{\underline \mu-b\leq y\leq \overline\mu-b
}\phi(y)=\phi(0)=1-{\rm e}^{-\varepsilon}>0.$$ Thus
$$\sum_{k=1}^\infty \CV\left(\left|\frac{ S_{n_k}-S_{n_{k-1}}}
{n_k-n_{k-1}}-b\right|\leq\varepsilon\right)\geq
\sum_{k=1}^\infty\CE\left[\phi\left(\frac{ \overline
S_{n_k-n_{k-1}}} {n_k-n_{k-1}}\right)\right]=\infty.$$ Note the fact
that the sequence of $\{S_{n_k}-S_{n_{k-1}}\}_{k\geq1}$ is
independent for all $k\ge 1$ . Using the second Borel-Cantelli
Lemma, we have
 $$
\liminf_{k\to\infty}\left|\frac{ S_{n_k}-S_{n_{k-1}} }
{n_k-n_{k-1}}-b\right|\leq \varepsilon,\quad \hbox{a.s.} \; \CV.
$$
But
\begin{equation}\label{ED}
\left|\frac{S_{n_k}}{n_k}-b\right|\leq \left|\frac{
S_{n_k}-S_{n_{k-1}}
}{n_k-n_{k-1}}-b\right|\cdot\frac{n_k-n_{k-1}}{n_k}+\left[\frac{|S_{n_{k-1}}|}{n_{k-1}}+|b|\right]\frac{{n_{k-1}}}{n_k}.
\end{equation}
Noting the following fact,
$$
\frac{n_k-n_{k-1}}{n_k}\to 1, \quad \frac{{n_{k-1}}}{n_k}\to 0,\;
\hbox{as} \; k\to \infty
$$
and
$$
\limsup_{n\to\infty} S_n/n\leq \overline \mu, \quad
\limsup_{n\to\infty} (-S_n)/n\leq -\underline \mu,
$$
which implies
$$
\limsup_{n\to\infty} |S_n|/n\leq \max\{|\overline \mu|, |\underline
\mu|\}<\infty.
$$
Hence, from inequality (16), for any $\varepsilon>0,$
$$
\liminf_{k\to\infty}\left|\frac{S_{n_k}}{n_k}-b\right|\leq\varepsilon,
\hbox{ a.s.}\; \CV.
$$ That is
$$\CV\left(\liminf_{n\to \infty} |S_{n}/{n}-b|\leq\varepsilon\right)=1.$$
Since $\varepsilon$ can be  arbitrarily small, we have

$$\CV\left(\liminf_{n\to \infty} |S_{n}/{n}-b|=0\right)=1.$$

Now we prove that $$\CV\left( \limsup\limits_{n\to\infty}S_n/n=
\overline \mu\right)=1,\quad \CV\left( \liminf\limits_{n\to\infty}
S_n/n= \underline \mu\right)=1.
$$ Indeed, from
$\CV\left(\liminf\limits_{n\rightarrow\infty}\left|\frac{S_{n}}{n}-b\right|=0\right)=1$,
$\forall b\in(\underline{\mu},\overline{\mu})$, we can obtain
$$\CV\left(\limsup\limits_{n\rightarrow\infty}\frac{S_n}{n}\geq b\right)=1.$$
Then
$$\CV\left(\limsup\limits_{n\rightarrow\infty}\frac{S_n}{n}\geq\overline{\mu}\right)=1.$$

On the other hand, we have, for each $P\in{\cal P}$,
$$P\left(\limsup\limits_{n\rightarrow\infty}\frac{S_n}{n}=
\overline{\mu}\right)=P\left(\limsup\limits_{n\rightarrow\infty}\frac{S_n}{n}\geq\overline{\mu}\right)
+P(\limsup\limits_{n\rightarrow\infty}\frac{S_n}{n}\leq\overline{\mu})-1.$$
But
$$P\left(\limsup\limits_{n\rightarrow\infty}\frac{S_n}{n}\leq\overline{\mu}\right)\geq
v\left(\underline{\mu}\leq\liminf\limits_{n\rightarrow\infty}\frac{S_n}{n}\leq\limsup\limits_{n\rightarrow\infty}\frac{S_n}{n}\leq\overline{\mu}\right)=1.$$
Thus
$$\CV\left(\limsup\limits_{n\rightarrow\infty}\frac{S_n}{n}=\overline{\mu}\right)=1.$$

In a similar manner, we can obtain
$$\CV\left(\liminf\limits_{n\rightarrow\infty}\frac{S_n}{n}=\underline{\mu}\right)=1.$$
The proof of Theorem 3.2 is complete.\\

\noindent{\large\bf 3.2 General weak law of large numbers for
capacities induced
by sublinear expectations}\\

\noindent{\bf THEOREM 3.3 (General weak law of large numbers for
capacities)}. Under the conditions of Theorem 3.1, then for any
$\varepsilon>0$,
\begin{equation}\label{ED}\lim\limits_{n\rightarrow\infty}v\left(\frac{1}{n}\sum_{i=1}^nX_i\in
(\underline\mu-\varepsilon,\overline\mu+\varepsilon)\right)=1.\end{equation}

\noindent{\bf PROOF.} To prove (17), we only need to prove that
$$\lim\limits_{n\rightarrow\infty}\CV\left(\frac{1}{n}\sum_{i=1}^nX_i\leq\underline\mu-\varepsilon\right)=0$$
and
$$\lim\limits_{n\rightarrow\infty}\CV\left(\frac{1}{n}\sum_{i=1}^nX_i\geq\underline\mu+\varepsilon\right)=0.$$
Now we prove
$$\lim\limits_{n\rightarrow\infty}\CV\left(\frac{1}{n}\sum_{i=1}^nX_i\leq\underline\mu-\varepsilon\right)=0.$$
For any $\delta<\varepsilon$, construct two functions $f$, $g$ such
that
$$f(x)=1\ \  \hbox{for}\ \ x\leq\underline\mu-\varepsilon-\delta,$$ $$f(x)=0\ \ \hbox{for}\ \
x\geq\underline\mu-\varepsilon,$$
$$f(x)=\frac{1}{\delta}(\underline\mu-\varepsilon-x)\ \ \hbox{for}\
\ \underline\mu-\varepsilon-\delta<x<\underline\mu-\varepsilon$$ and
$$g(x)=1\ \ \hbox{for}\ \ x\leq\underline\mu-\varepsilon,$$ $$g(x)=0\ \ \hbox{for}\ \ x\geq
\underline\mu-\varepsilon+\delta,$$
$$g(x)=\frac{1}{\delta}(\underline\mu-\varepsilon-x)+1\ \ \hbox{for}\ \
\underline\mu-\varepsilon<x<\underline\mu-\varepsilon+\delta.$$
Obviously, $f$ and $g\in C_{b,Lip}(\CR)$. By Theorem 3.1, we have
$$\lim\limits_{n\rightarrow\infty}\CE\left[f\left(\frac{1}{{n}}\sum\limits_{i=1}^nX_i\right)\right]
=\sup\limits_{\underline \mu\leq y\leq\overline\mu}f(y)=0,$$
$$\lim\limits_{n\rightarrow\infty}\CE\left[g\left(\frac{1}{{n}}\sum\limits_{i=1}^nX_i\right)\right]
=\sup\limits_{\underline \mu\leq y\leq\overline\mu}g(y)=0.$$ So
\begin{equation}\label{ED}0=\sup\limits_{\underline \mu\leq y\leq\overline\mu}f(y)\leq\liminf\limits_{n\rightarrow\infty}
\CV\left(\frac{1}{n}\sum_{i=1}^nX_i\leq\underline\mu-\varepsilon\right)\leq\limsup\limits_{n\rightarrow\infty}
\CV\left(\frac{1}{n}\sum_{i=1}^nX_i\leq\underline\mu-\varepsilon\right)\leq
\sup\limits_{\underline \mu\leq
y\leq\overline\mu}g(y)=0.\end{equation} Hence
$$\lim\limits_{n\rightarrow\infty}\CV\left(\frac{1}{n}\sum_{i=1}^nX_i\leq\underline\mu-\varepsilon\right)=0.$$

In a similar manner, we can obtain
$$\lim\limits_{n\rightarrow\infty}\CV\left(\frac{1}{n}\sum_{i=1}^nX_i\geq\underline\mu+\varepsilon\right)=0.$$
The proof of Theorem 3.3 is complete.


\begin{thebibliography}{99}

\bibitem{Chen} Chen, Z. J., Strong laws of large numbers for capacities, arXiv:1006.0749v1 [math.PR] 3 Jun 2010.


\bibitem{Lei} Li, M. and Shi, Y. F., A general central limit theorem under
sublinear expectations, Science China Mathematics, {\bf 53}, 2010,
1989-1994.

\bibitem{Maccheroni} Maccheroni, F. and Marinacci, M.,  A strong law of large number for capacities, Ann.
Proba.,  {\bf 33}, 2005, 1171-1178.


\bibitem{Marinacci} Marinacci, M.,  Limit laws for non-additive
probabilities and their frequentist interpretation, J. Econom.
Theory, {\bf 84}, 1999, 145-195.


\bibitem{Peng1} Peng, S. G., G-expectation, G-Brownian motion and related stochastic
calculus of Ito type, in: F.E. Benth, et al. (Eds.),  Proceedings
of the Second Abel Symposium, 2005, Springer-Verlag, 2006,
541-567.

\bibitem{Peng2} Peng, S. G., Law of large number and central limit theorem
under nonlinear expectations, arXiv: math. PR/0702358vl 13 Feb 2007.

\bibitem{Peng3} Peng, S. G., Multi-dimensional G-Brownian motion and related
stochastic calculus under G-expectation, Stochastic Process.
Appl., {\bf 118(12)}, 2008, 2223-2253.

\bibitem{Peng4} Peng, S. G., A new central limit theorem under sublinear expectations,
arXiv:0803.2656v1 [math.PR] 18 Mar. 2008.

\bibitem{Peng5} Peng, S. G., Survey on normal distributions, central limit theorem,
 Brownian motion and the related stochastic calculus under sublinear
 expectations, Science in China Series A-Mathematics, {\bf 52(7)},
 2009, 1391-1411.

\bibitem{Wang}Wang, L. H., On the regularity of fully nonlinear parabolic equations II, Comm. Pure Appl. Math., {\bf
45}, 1992, 141-178.
\end{thebibliography}
\end{document}